\documentclass[12pt]{amsart}
\usepackage{graphicx}
\usepackage{amssymb}
\usepackage{epstopdf}
\newcommand{\rat}{\mathfrak{rat}}
\newcommand{\U}{u}
\renewcommand{\sl}{\mathfrak{sl}_2}
\newcommand{\A}{\mathcal A}

\newcommand{\Z}{\mathbb Z}
\newcommand{\two}{\A^{2h}_{n+1}(*)}

\newtheorem{theorem}{Theorem}

\newtheorem{prop}{Proposition}

\begin{document}
    \baselineskip=16pt - 

\title{Gropes and the rational lift of the Kontsevich integral}
\author{James Conant}
\address{Dept. of Mathematics\\
                 University of Tennessee\\
                 Knoxville, TN 37996-1300}

\email{jconant@math.utk.edu}
\thanks{Supported by NSF grant DMS 0305012}
\subjclass{57M27}

\begin{abstract}
In this note, we calculate the leading term of
the rational lift of the Kontsevich integral, $Z^{\rat}$, introduced by Garoufalidis and Kricker, on the boundary of an embedded grope of class $2n$. We observe that it lies in the subspace spanned by connected diagrams of Euler degree $2n-2$ and with a bead $t-1$ on a single edge. This places severe algebraic restrictions on the sort of knots that can bound gropes, and in particular
implies the two main results of the author's thesis \cite{c}, at least over the rationals.  
\end{abstract}

\maketitle

\section{Introduction}
\begin{figure}
\includegraphics[height=1.7in]{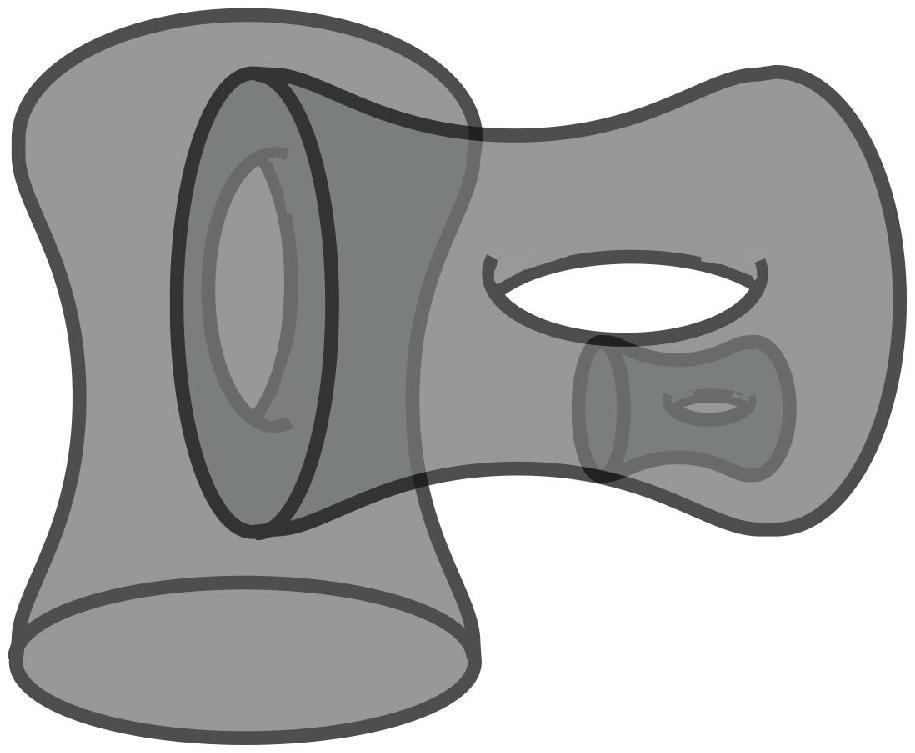}
\includegraphics[height=1.7in]{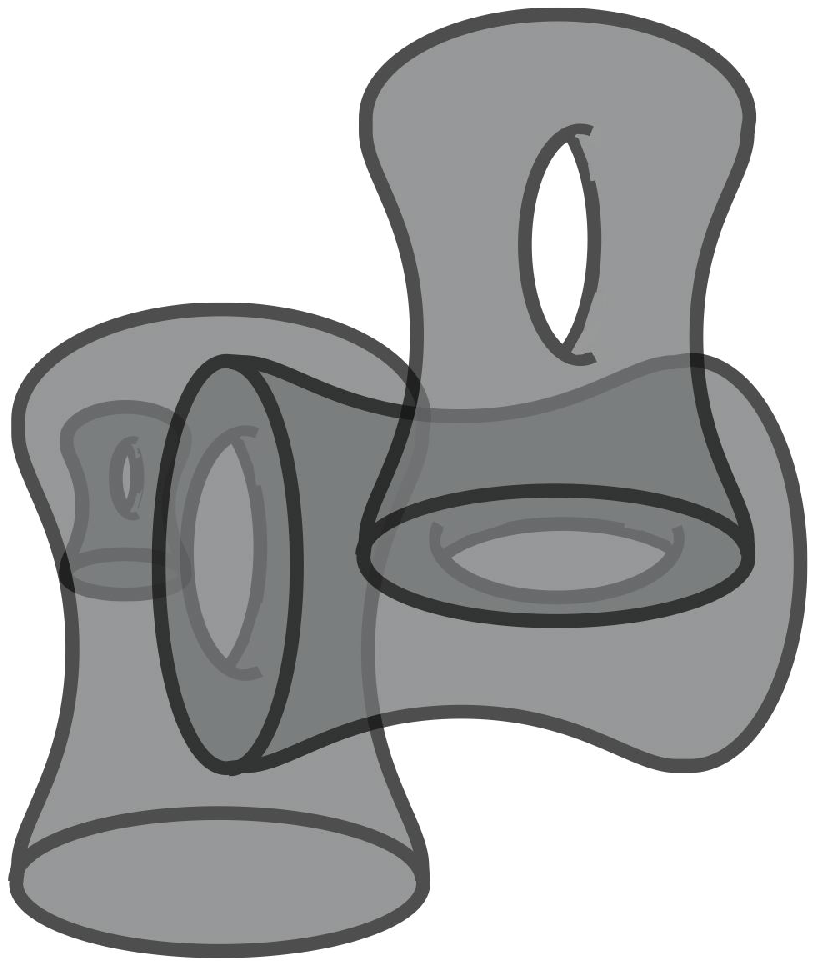}
\includegraphics[height=1.7in]{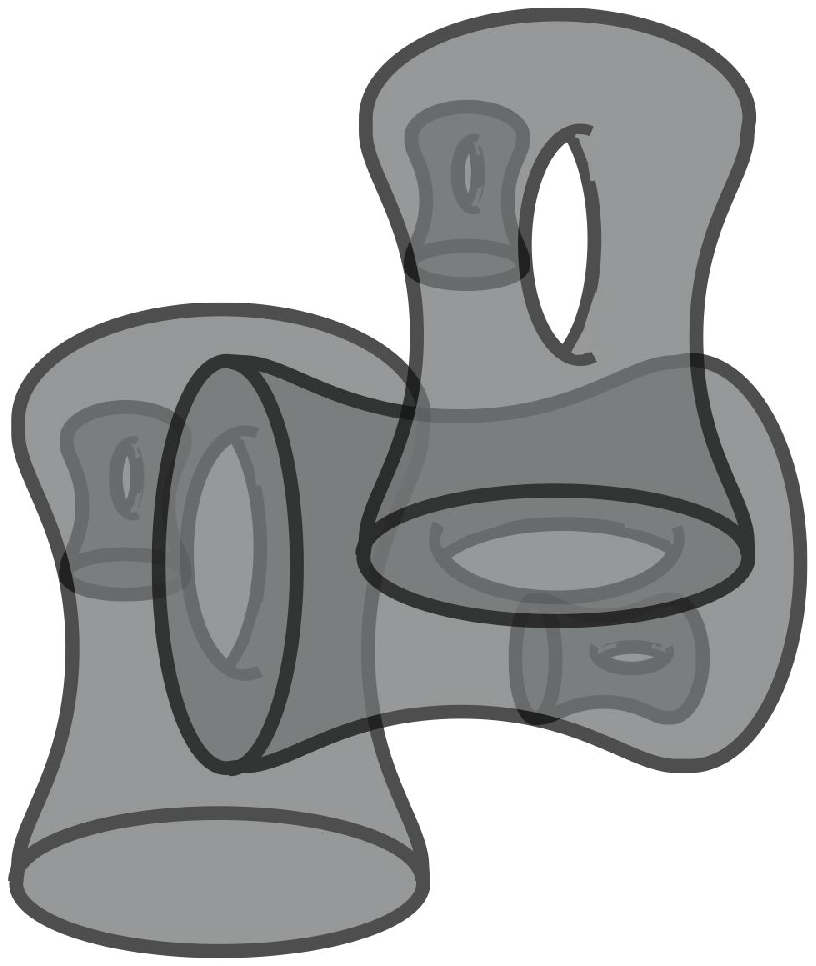}
\caption{Gropes of class $4$,  $5$ and $7$.}\label{fig1}
\end{figure}

In \cite{gk}, Garoufalidis and Kricker introduced a powerful invariant $Z^{\rat}$ which is a lift of the Kontsevich integral to a space spanned by trivalent graphs whose edges are colored by elements of $\Z[t,t^{-1}]$, at least when the Alexander polynomial vanishes. In this note we calculate how $Z^{\rat}$ behaves when evaluated on the boundary of a grope. A grope (of class $c$) is a specific $2$-complex 
that topologically models an element of the $c$th term of the lower central series of a group. It is formed by gluing together surfaces called \emph{stages.}
Some examples are pictured in Figure~\ref{fig1}.
See \cite{ct1} and \cite{teich} for details.

\begin{theorem}\label{rational}
 Let $g$ be an embedded grope of class $2n$ $(n\geq 1)$ in $S^3$ which has one boundary component, is of class $2n$, and has surface stages of genus one. Let $\U$ denote the unknot.
 Then
\begin{itemize}
 \item[a)] $Z^{\rat}(\partial g-\U)$ vanishes in Euler degree below $2n-2$.
 \item[b)] $Z_{2n-2}^{\rat}(\partial g-\U)$ is a sum of connected graphs where every edge but one is colored by $1\in \Z[t,t^{-1}]$, and the other edge is colored by $t-1$. 
\end{itemize}
\end{theorem} 

Note that the knots in this theorem automatically have trivial Alexander polynomial since they bound gropes of class exceeding three.

This will imply that the Kontsevich integral of $\partial g-\U$ vanishes in degree $\leq n$ and that
$$Z_{n+1}(\partial g-\U)\in\two$$
where $\two$ denotes the subspace of $\A_{n+1}(*)$ spanned by connected graphs with two hairs on a single edge. Moreover we will show every element of $\two$ is realized for some $g$.
Since one can easily show that $\two\neq 0$ using an $\sl$ weight system, we achieve the following results, at least when the ambient manifold is $S^3$. (These were both proven in \cite{c} with a lot more effort.)

\begin{theorem}\label{vassiliev}
Let $g$ be an embedded genus one grope of class $2n$ in a homology $3$-sphere ($n\geq 1$).
\begin{itemize}
\item[a)] Rational-valued Goussarov-Vassiliev invariants of degree $\leq n$ do not distinguish $\partial g$
from the unknot $\U$. 
\item[b)] There exists a Goussarov-Vassiliev invariant, $v$, of degree $n+1$ and a grope $g$ of class $2n$ such that $v(\partial g)\neq v (\U)$.
 \end{itemize}
\end{theorem}

The case of arbitrary homology spheres follows easily from the $S^3$ case, as we show at the end of the paper.

We conclude by emphasizing the fact that this note deals with knots that bound gropes, which is a more restricted class than those cobounding a grope with the unknot, as studied in \cite{ct1}. 

\section{Evaluating $Z^{\rat}$ on the boundary of a grope}
\subsection{The invariant $Z^{\rat}$}
We briefly review some properties of $Z^{\rat}$ here. This is distilled from \cite{gt}, which is further distilled from \cite{gk}. We restrict to the case of knots with trivial Alexander polynomial, since knots which bound gropes of class $\geq 3$ have that property.

Let $\Lambda=\Z[t,t^{-1}]$. Then a $\Lambda$-colored graph is a trivalent graph together with a map
$\operatorname{Edges}\to \Lambda$. (A coloring of an edge is sometimes called a bead.)
Then $\A(\Lambda)$ is a vector space generated by $\Lambda$-colored graphs, modulo certain relations. This space is graded by the Euler degree which is the number of vertices of the graph.

The so-called hair map
$$\operatorname{Hair}\colon\A(\Lambda)\to \A(*)$$
is defined by expanding the variable $t$ coloring an edge $e$ into an infinite series $\operatorname{exp}(h)\cdot e$,
where $h\cdot e$ is by definition the addition of a single hair.

Then $\operatorname{Hair}\circ Z^{\rat}=Z$, indicating that $Z^{\rat}$ is indeed a lift of $Z$.

Letting $Z_{2n}^{\rat}$ be the degree $2n$ part of $Z^{\rat}$, we will use the fact that
$Z_{2n}^{\rat}$ is a universal $\mathbb Q$-valued finite type invariant of null degree $2n$. That is, it is universal with respect to null clasper surgeries, which are  clasper surgeries whose leaves link the knot trivially.

So $Z^{\rat}$ vanishes on alternating sums of surgeries of  null claspers of total degree $2n-1$, and its value on alternating sums of null-claspers of total degree $2n$ is given by the so-called \emph{complete contraction}, defined by gluing the leaves together using the equivariant linking number, and coloring an edge by a $t$ or $t^{-1}$ every time it passes through a fixed Seifert surface. 

The equivariant linking number of two leaves which are meridians to dual bands of a Seifert surface is $t-1$. More generally, suppose that one leaf is a meridian to a band of a Seifert surface and another leaf
links only the dual band with linking number $\ell$. Then the equivariant linking number of these two leaves is $\ell\cdot(t-1)$. 
On the other hand, the equivariant linking number of leaves which are null homotopic in the complement of a Seifert surface is just the standard linking number. 

\subsection{Proof of Theorem \ref{rational}}

Suppose a knot bounds an embedded grope of class $2n$ in $S^3$, where all the surface stages are of genus one. 
Then the knot can be represented as a rooted tree clasper surgery, $T$, in the complement of the unknot which forms a meridian to the root leaf \cite{ct1}. 
Note that the other leaves can be embedded arbitrarily in the complement of the unknot's spanning disk.
 Break $T$ into a union of a $Y$ at its root and a clasper $T^\prime$ which sits on the unknot $\U_Y$.
Then $T^\prime$ is a null clasper of Euler degree $2n-2$. The complete contraction $\langle T^\prime \rangle$ is defined by thinking of $T^\prime$ as a unitrivalent tree and gluing up its leaves using the equivariant linking form.
One of the two dual bands of $\U_Y$ has exactly one leaf of $T^\prime$ as a meridian. (Because the grope is of class exceeding $2$.)  The equivariant linking number of this leaf with any leaf that links the dual band is an integral multiple of $t-1$. 
 For every other pair of leaves the equivariant linking numbers are integers.
Thus $\langle T^\prime\rangle$ is a sum of trivalent graphs of Euler degree $2n-2$ which have a bead $t-1$ on a single edge. By the universality of $Z^{\rat}$, we have that $Z^{\rat}_{<2n-2}(\U_T-\U)=0$, and that
$$Z^{\rat}_{2n-2}(\U_T-\U)=\langle T^\prime\rangle.$$

\subsection{Two propositions}

\begin{prop}
The Kontsevich integral of $\partial g-\U$ vanishes in degree $\leq n$ and 
$Z_{n+1}(\partial g-\U)\in\two.$
Moreover, every element of $\two$ is equal to $Z_{n+1}(\partial g-\U)$ for some $g$.
\end{prop}
\begin{proof}
To calculate $Z$, we need to apply the hair map:
\begin{align*}
\operatorname{Hair}\langle T^\prime\rangle&=  h\cdot \langle T^\prime\rangle+(1/2)h^2\cdot \langle T^\prime \rangle+\cdots\\
      						      &=  (1/2)h^2\cdot \langle T^\prime\rangle+\cdots
\end{align*}
This follows since any trivalent graph with a single univalent vertex is trivial in $\A(*)$.

Note that any terms of $Z^{\rat}$ of Euler degree exceeding $2n-2$ must have at least two hairs and so be of Vassiliev degree exceeding $n+1$. 
Thus there are no terms with Vassiliev degree less than or equal to $n$, and
 $Z_{n+1}(\U_T-\U)=(1/2)h^2\cdot \langle T^\prime \rangle$.

Observe that $(1/2)h^2\cdot \langle T^\prime\rangle\in\two$ as claimed. To see that every element of $\two$ is realized, take a generator $D\in\two$. Remove one of the two hairs and embed $D$ as a rooted clasper on the unknot. Breaking some edges of $D$ into hopf pairs of leaves we get a rooted clasper with $2n$ leaves as desired.
\end{proof}

\begin{prop}
For $n\geq 1$ we have $\two\neq 0$.
\end{prop}
\begin{proof}
This can be detected using an $\sl$ weight system \cite[Appendix]{t}. This weight system satisfies two relations. Any closed loop gives a multiplicative factor of $3$, whereas an edge can be expanded into a difference of two terms, where the edge is removed and replaced either by two parallel edges, or by two edges which cross.

Consider the diagram 

\vspace{.2in}
\centerline {
\includegraphics[width=1.5in]{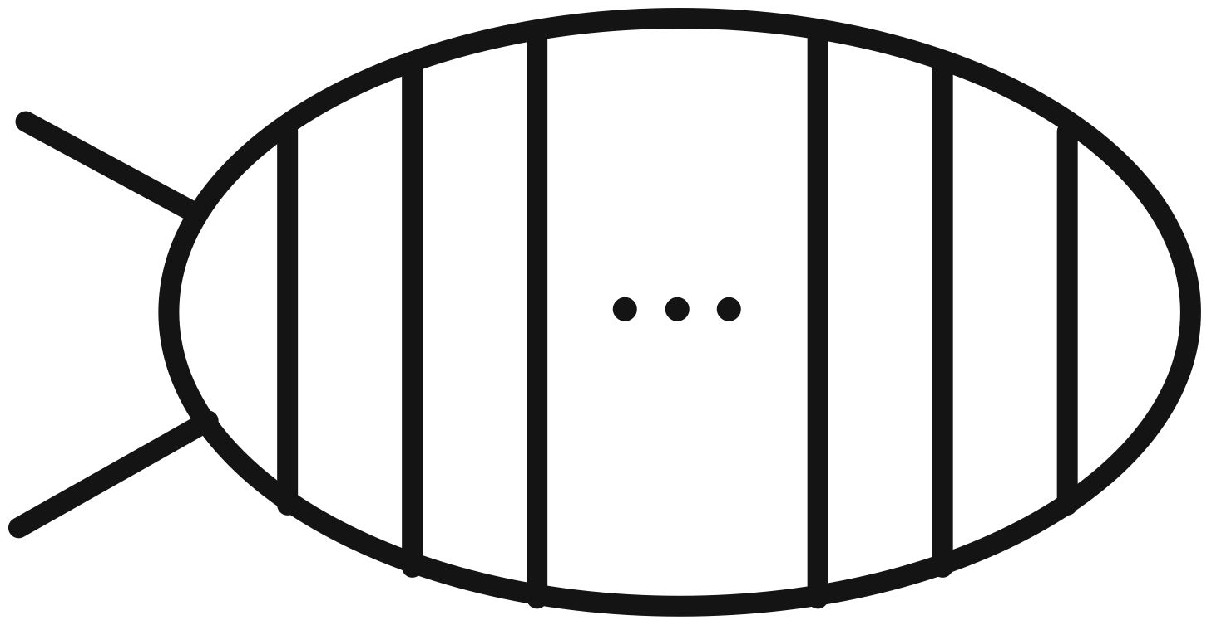}
}
\vspace{.2in}
\noindent which consists of an ellipse with $n-1$ parallel vertical edges and two hairs. We show it is nonzero by showing that the diagram formed by joining the two
univalent vertices is nonzero. The resulting diagram is an ellipse with $n$ parallel edges added.
Applying the $\sl$ relations to these $n$ edges, starting from the right, we get down to $2^{n}\cdot\mathbf O$, where $\mathbf O$ is a single circle. Hence our diagram evaluates to $2^n\cdot 3$.
\end{proof}

\subsection{Proof of Theorem \ref{vassiliev}}

We have already established Theorem 1 for knots in $S^3$. Any homology sphere, $\Sigma$, can be obtained from $S^3$ by surgery on a disjoint union of $Y$ claspers. By general position, we may assume that the clasper $T$ corresponding to our grope is disjoint from these claspers. 
Letting $T^\prime$ denote the clasper formed by surgering the unknot along the $Y$ containing the root, we see that modulo Euler degree $2n-1$ we can pull $T^\prime$ into a ball. Then $Z^{\rat}$ decomposes as a product of $Z^{\rat}(\Sigma)$ and $Z^{\rat}(\U_T)$, which easily implies Theorem 1.

 \end{document}